\documentclass[amssymb,aps,showpacs,showkeys,reprint]{revtex4-1}

\usepackage{graphicx}
\usepackage{amsmath}
\usepackage{color}
\definecolor{brown}{rgb}{0.8,0.6,0.3}
\definecolor{dgreen}{rgb}{0.2,0.4,0.3}

\usepackage{fancyvrb}
\usepackage[pdfauthor={Richard J. Mathar},colorlinks=true,citecolor=brown,urlcolor=blue,linkcolor=red,pdfkeywords={Simple Cubic Grid, Subgrid, Tessellation},pdfpagelayout=OneColumn]{hyperref}

\begin{document}

\title{Hierarchical Subdivision of the Simple Cubic Lattice}

\author{Richard J. Mathar}
\homepage{http://www.mpia.de/~mathar}
\email{mathar@mpia.de}
\affiliation{Hoeschstr. 7, 52372 Kreuzau, Germany}

\date{\today}

\begin{abstract}
The simple cubic lattice defines a set of points at regular distances. The 
volume of the Voronoi
cells around each point may serve as a weight
for integration over the entire space. We add interstitial
points to this grid according to the rule that these have maximum distance
to the existing points, or by the equivalent rule that they are placed at the vertices of
the Voronoi cells of the existing lattice.
Choices
of that kind appear in numerical sampling with maximum independent data points
if the numerical expense of computing values represented at the lattice points is high.
The volumes and shapes of the
Voronoi cells of these enriched/supersampled lattices are discussed in detail
while this insertion is recursively
executed three times in succession.
\end{abstract}

\maketitle
\section{Theme}
Some numerical schemes sampling space with finite elements for interpolation or integration
may employ mesh refinement techniques which start from
a coarse grid of mesh points and build finer grids with denser point sets
at latter stages. The requirements of (i) efficient use of the values samples
at previous stages, and (ii) maximum independent information contributed by the
refined stages will be embodied in this manuscript by the following recipe
and guideline:
Place points of the new, refined stages at interstitial positions of the previous grid such
that they have maximum Euclidean distance to the points of the previous grid.

Other criteria based on Fourier amplitudes are possible and have been applied to 
symmetry-adapted sampling of the reciprocal lattice \cite{CunninghamPRB10,PackPRB16,EntezariIEEEvis04}.

This idea will be worked out by starting from the square or hexagonal
grid in two
dimensions and as warm-up, then the simple cubic lattice in three dimensions.
These grids are the basic ones because coding the lattice points
in Cartesian coordinates is as simple as running with integer
lattice coordinates independently through the base vectors of their lattices.

\section{Two-dimensional Templates} \label{sec.2d}
\subsection{Square Grid}
The simplest coverage of the plane by a point grid is the square grid with
lattice points placed at
\begin{equation}
\mathbf{p} = i \mathbf{e}^{(0)}_1 + j\mathbf{e}^{(0)}_2,
\end{equation}
where the two unit vectors $\mathbf{e}_{1,2}^{(0)}$ point from
the Cartesian coordinates $(0,0)$ at
the origin to $(a,0)$ and $(0,a)$. $a$ is the lattice constant,
and $i$ and $j$ run through all integers of both signs.

There are two evident choices for the unit cell:
\begin{itemize}
\item
The \emph{primitive}
unit cell is a square
and stretches from one lattice point to the two nearest
neighbors at distance $a$ to the right and up, and along
diagonals of the square to the second nearest neighbors at distance $\sqrt{2}a$.
\item
The \emph{Wigner-Seitz} (WS) (or Voronoi) unit cell uses the Brillouin zone construction
of solid state physics. Between each
pair of points in the grid, a line is drawn that cuts mid-way orthogonally through the
finite line that connects the two points, and which divides space
into
two half-spaces. The closest polygon around a point
build from the finite pieces of these lines of separation encompasses
the unit cell. It is the convex polytope defined by intersection
of the half-spaces.  For the square grid, this unit cell is also a square of edge
length $a$, but with a lattice point in the square center.
\end{itemize}
These unit cells have area $a^2$.

The obvious choice for an additional set of points towards a denser mesh
is to add one new point at the center of each primitive unit cell. This meets the 
requirement of maximum distance formulated above, because these points are
$a/\sqrt{2}$ away from any of the points of the original mesh. So around
these we can draw the largest circles that contain no point of the original
mesh.

These points could also be found by putting them at the
vertices of the WS unit cells with maximum degree (the degree being the
number of WS unit cells that meet there). By construction, these vertices
are as far away from as many as possible points of the original grid. The
reasoning with the Voronoi cells is more satisfactory from a conceptional
point of view: finding a point in the plane with largest free circumcircle
tastes like a minimization problem that needs numerical treatment.
On the other hand,
given a point $\mathbf{o}$ in the plane with surface  normal $\mathbf{n}$,
$\mathbf{o}$ has a unique representation
$\mathbf{o}=\alpha\mathbf{n}+\mathbf{r}$ with $\mathbf{r}\perp\mathbf{n}$,
$\alpha=\mathbf{o}\cdot \mathbf{n}$, and with normal form  $Ax+By+Cz+D=0$, $D=-\alpha$
of the plane, where $A$, $B$, $C$ are the Cartesian components of $\mathbf{n}$.
Finding the common vertex of an intersection of three planes is therefore
equivalent to solving a $3\times 3$ inhomogeneous linear system of equations
for three unknowns $x,y$ and $z$. (The equivalent statements appear in the 2-dimensional problem.)
With some knowledge of which three lattice points nearby a pivotal lattice point
define a vertex of the Voronoi polytope, the explicit determination of the Cartesian
components of the vertex is therefore easy.

The additional grid points are obviously a copy of the original grid
points translated by half a diagonal. They are located at positions
\begin{equation}
(i+\frac{1}{2})\mathbf{e}^{(0)}_1+(j+\frac{1}{2})\mathbf{e}^{(0)}_2
\end{equation}
with integer $i$ and $j$. A pleasant observation is that the union
of the old and new grid points \emph{again} represents a simple square
grid---with smaller lattice constant $a/\sqrt{2}$ and rotated by $45^\circ$
relative to the original one. This could be noted by defining the unit vectors
\begin{equation}
\mathbf{e}^{(1)}_1\equiv \frac{1}{2}(\mathbf{e}^{(0)}_1+\mathbf{e}^{(0)}_2),\quad
\mathbf{e}^{(1)}_2\equiv \frac{1}{2}(-\mathbf{e}^{(0)}_1+\mathbf{e}^{(0)}_2)
\end{equation}
and accessing the union of these lattice points of two levels by
\begin{equation}
i\mathbf{e}^{(1)}_1+j\mathbf{e}^{(1)}_2
.
\end{equation}
This leads to a simple recursive mesh refinement: iteratively
the lattice constant is divided by $\sqrt 2$, and the
two orthogonal unit vectors are rotated alternatingly by $45^\circ$ to the left
and to the right.

\subsection{Triangular (Hexagonal) Grid}
The fundamental areal element of the triangular grid is the isosceles triangle
with side length $a$, and two unit vectors with Cartesian coordinates
\begin{equation}
\mathbf{e}^{(0)}_1= (a,0),\quad \mathbf{e}^{(0)}_2 = (-a/2,\sqrt{3}a/2)
\end{equation}
with an angle of $120^\circ$ between.
The unit cell has an area of $\sqrt{3}a^2/4$. The primitive unit cell
is a kite shaped quadrangle formed by two of these triangular elements glued by one side.
The WS unit cell is a hexagon centered at a lattice point.

The refinement of this triangular lattice is again obvious. By any of the two
methods proposed above (maximum free distance or vertices of the WS unit cell)
the additional points are placed in the middle (mid-point of the circumcircle)
of each the two triangles of the primitive unit cell, called the $K$ point
in solid state plane groups \cite{TerzibaschPSS133}.
The first has Cartesian
coordinates $(a/2,a/(2\sqrt{3}))$, the second $(0,a/\sqrt{3})$.

The union of the grid points
of the original lattice and the additional points at the interstitial locations
establishes another triangular lattice, rotated by $30^\circ$
from the original lattice. The lattice constant is shrunk
to $a/\sqrt 3$, the nearest neighbor distance in the refined grid.
So the area of the unit cell of the refined grid has shrunk by a factor $1/3$
compared to the area of the original cell. The factor of three is essentially
indicating that the density of the mesh points has triplicated because
we added \emph{two} grid points into each unit cell of the original mesh,
which contained \emph{one} grid point.

Similar to the finding with the square grid, this refinement preserves the grid structure, and recursive
refinement is therefore a rather easy task from a programmer's point of view. 
At each step, the previous unit vectors are shrunk by a factor $1/\sqrt 3$
and rotated by $30^\circ$.

\section{Three-dimensional Grid}
\subsection{Level 0: Simple Cubic}

The starting point for an unbiased sampling
in three dimensions is the Simple Cubic (SC) lattice with lattice constant $a$,
unit vectors $\mathbf{e}^{(0)}_1=(a,0,0)$
$\mathbf{e}^{(0)}_2=(0,a,0)$
$\mathbf{e}^{(0)}_3=(0,0,a)$.
The volume of the unit cell is
\begin{equation}
V^{(0)}=a^3.
\label{eq.scv}
\end{equation}
Points are located at positions
\begin{equation}
i \mathbf{e}^{(0)}_1
+ j \mathbf{e}^{(0)}_2
+ k \mathbf{e}^{(0)}_3
\label{eq.scp}
\end{equation}
with integer coordinates $i$, $j$ and $k$.

Counting neighbors in shells of common distance to the grid points is
a useful digital signature of the grid structure. For the SC
lattice this statistics starts as in Table \ref{tab.nl0} for the
smallest distances. Top to bottom,
the 6 nearest neighbors are in the directions of the unit
vectors $\mathbf{e}^{(0)}_i$ and their opposites,
the 12 second nearest neighbors are found at positions along the face diagonals, and
the 8 third nearest neighbors are
in the directions of the space diagonals \cite[A005875]{EIS}.
\begin{table}
\caption{Frequencies of distances to neighbors in the SC lattice.}
\begin{tabular}{r|r}
\hline
count & squared distance\\
\hline
6 & $a^2$ \\
12 & $2a^2$ \\
8 & $3a^2$ \\
6 & $4a^2$ \\
24 & $5a^2$ \\
24 & $6a^2$ \\
\hline
\end{tabular}
\label{tab.nl0}
\end{table}

\subsection{Level 1: Body-centered Cubic} \label{sec.bcc}

The obvious first refinement places interstitial points half way
along the space diagonal for another copy with points at
\begin{equation}
(i+\frac{1}{2})
\mathbf{e}^{(0)}_1
+
(j+\frac{1}{2})
\mathbf{e}^{(0)}_2
+
(k+\frac{1}{2})
\mathbf{e}^{(0)}_3
\label{eq.bccp}
\end{equation}
with integer triples $i$, $j$ and $k$.

The major difference in comparison with the two-dimensional examples
is that the union of these grid points does not define another
simple cubic lattice but a body-centered cubic (BCC) lattice.
The WS unit cell of this lattice is characterized by six
squares in the directions of the six next nearest neighbors
of the original unit cell plus eight hexagons in the directions
of the eight new lattice points in the centers of the eight
primitive cells with common vertex at $(0,0)$.
In Brillouin zones of the corresponding space group,
the mid points of the squares are labeled $X$
and the mid points of the hexagons are labeled $L$, and
the vertices of the squares and hexagons are labeled $W$
\cite{Bouckaert,HerringPR52,ElliottPR96}.

The statistics of neighbors around the lattice points
(at both levels) is indicated in Table \ref{tab.nl1}.
\begin{table}
\caption{Statistics of distances for points in the BCC lattice. \cite[A004013]{EIS}}
\begin{tabular}{r|r}
\hline
count & squared distance \\
\hline
8 & $\frac34 a^2$\\
6 & $a^2$ \\
12 & $2a^2$ \\
24 & $\frac{11}{4}a^2$ \\
8 & $3a^2$ \\
6 & $4a^2$ \\
24 & $\frac{19}{4}a^2$ \\
24 & $5a^2$ \\
24 & $6a^2$\\
\hline
\end{tabular}
\label{tab.nl1}
\end{table}
The frequencies of neighbors around any of the new
points at the positions (\ref{eq.bccp})
is the same as for the points of the
zeroth level of refinement (\ref{eq.scp}).

\begin{figure}
\includegraphics[width=0.99\columnwidth]{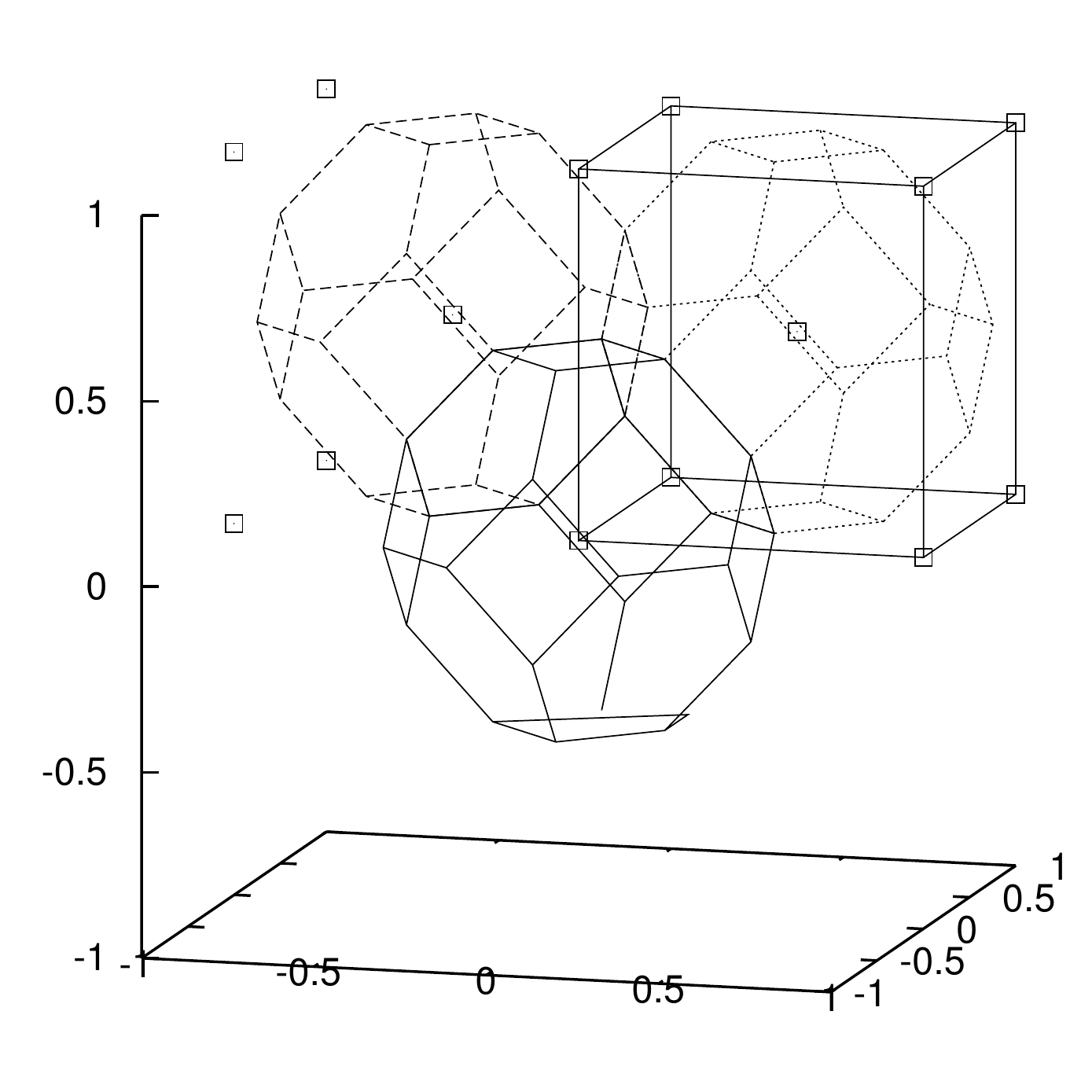}
\caption{Points on a BCC lattice (dotted squares). Three of them are surrounded
by the truncated octahedra of the Wigner-Seitz cell.}
\label{fig.WSbcc}
\end{figure}

\begin{figure}
\includegraphics[width=0.99\columnwidth]{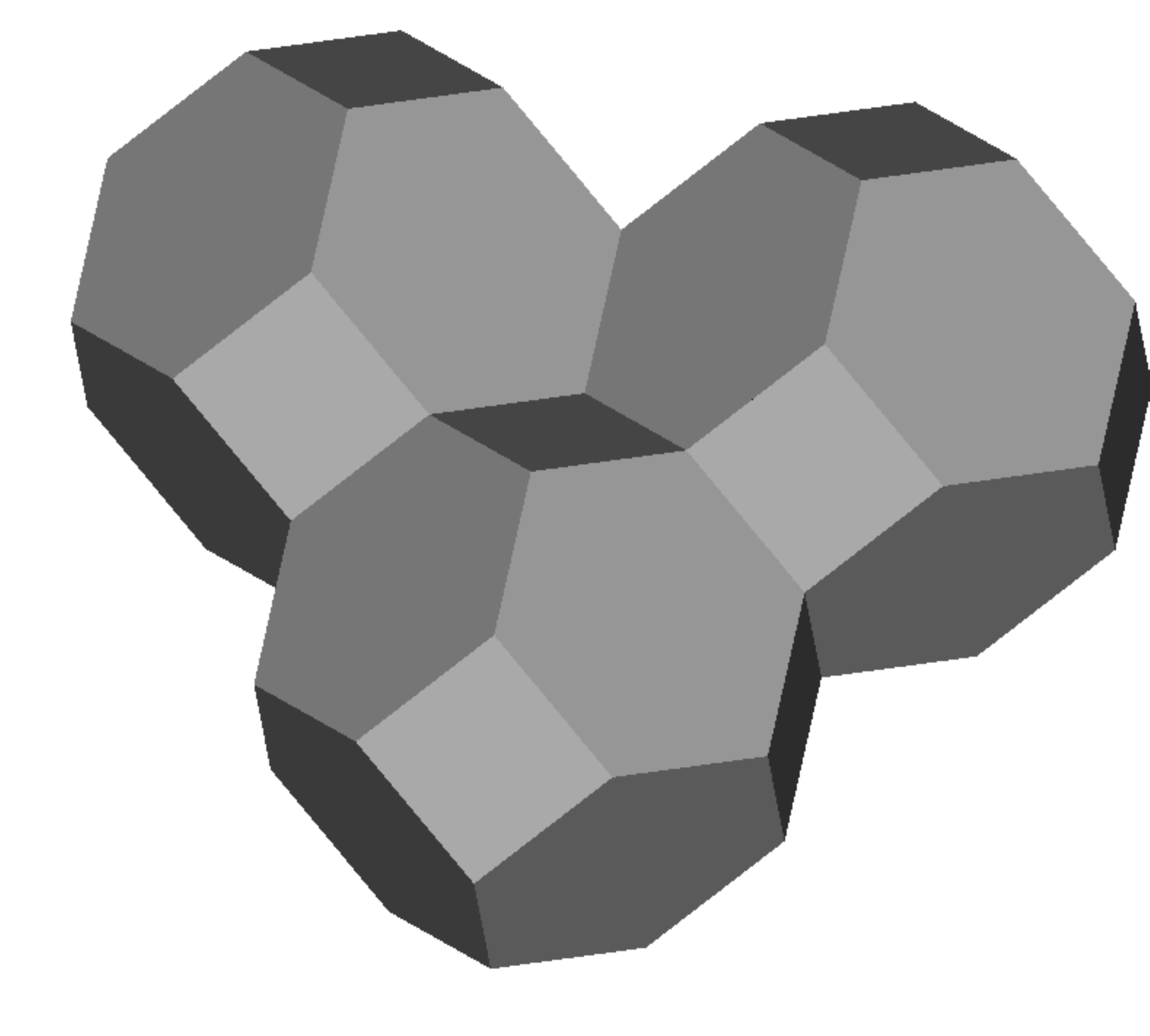}
\caption{WS cells of the BCC lattice as in Figure \ref{fig.WSbcc}, but
with back surfaces hidden.}
\label{fig.bcc}
\end{figure}

The weights associated with numerical integration are the volumes
of the Voronoi cells around each lattice point.
The two Voronoi cells
around the points at level 0---enumerated by (\ref{eq.scp})---
and around the points at level 1---enumerated by (\ref{eq.bccp})---have
the same shape, illustrated in Figure \ref{fig.WSbcc} and \ref{fig.bcc}.
The polyhedron of the Wigner-Seitz cell
of that BCC lattice is a Truncated Octahedron with six quadratic and eight regular hexagonal faces.
Figure \ref{fig.WSbcc} shows three of these unit cells. The two cells at the back
share a common square face; the one at the front is attached to the  one right at the back
by a hexagon.
The unit box is the  primitive unit cell of the SC lattice.
The volume of each Truncated Octahedron is half of the primitive unit cell,
\begin{equation}
V_\Gamma^{(1)}=a^3/2.
\label{eq.bccv}
\end{equation}
The upper index indicates that the volume is defined at level 1, after
one refinement, and the lower label is the standard symbol for
the type of symmetry of the center point of the cell.

\subsection{Level 2: body-centered with W or X}
\subsubsection{Level 2: body-centered with W} \label{sec.bcW}

The next refinement step adds points at $W$ positions
(the name of the points at the 24 vertices of the truncated octahedra
in the associated Brillouin zone of the
face-centered cubic lattice).

The coordinate triples of $W$ on the left, bottom and front faces of the
primitive unit cell of the SC are
\begin{multline}
i\mathbf{e}_1^{(0)}
+(j+\frac{1}{4}) \mathbf{e}_2^{(0)}
+(k+\frac{1}{2}) \mathbf{e}_3^{(0)},\\
i\mathbf{e}_1^{(0)}
+(j+\frac{3}{4}) \mathbf{e}_2^{(0)}
+(k+\frac{1}{2}) \mathbf{e}_3^{(0)},\\
i\mathbf{e}_1^{(0)}
+(j+\frac{1}{2}) \mathbf{e}_2^{(0)}
+(k+\frac{1}{4}) \mathbf{e}_3^{(0)},\\
i\mathbf{e}_1^{(0)}
+(j+\frac{1}{2}) \mathbf{e}_2^{(0)}
+(k+\frac{3}{4}) \mathbf{e}_3^{(0)},
\\
(i+\frac{1}{4})\mathbf{e}_1^{(0)}
+(j+\frac{1}{2}) \mathbf{e}_2^{(0)}
+k \mathbf{e}_3^{(0)},\\
(i+\frac{3}{4})\mathbf{e}_1^{(0)}
+(j+\frac{1}{2}) \mathbf{e}_2^{(0)}
+k \mathbf{e}_3^{(0)},\\
(i+\frac{1}{2})\mathbf{e}_1^{(0)}
+(j+\frac{1}{4}) \mathbf{e}_2^{(0)}
+k \mathbf{e}_3^{(0)},\\
(i+\frac{1}{2})\mathbf{e}_1^{(0)}
+(j+\frac{3}{4}) \mathbf{e}_2^{(0)}
+k \mathbf{e}_3^{(0)},
\\
(i+\frac{1}{4})\mathbf{e}_1^{(0)}
+j\mathbf{e}_2^{(0)}
+(k+\frac{1}{2}) \mathbf{e}_3^{(0)},\\
(i+\frac{3}{4})\mathbf{e}_1^{(0)}
+j\mathbf{e}_2^{(0)}
+(k+\frac{1}{2}) \mathbf{e}_3^{(0)},\\
(i+\frac{1}{2})\mathbf{e}_1^{(0)}
+j\mathbf{e}_2^{(0)}
+(k+\frac{1}{4}) \mathbf{e}_3^{(0)},\\
(i+\frac{1}{2})\mathbf{e}_1^{(0)}
+j\mathbf{e}_2^{(0)}
+(k+\frac{3}{4}) \mathbf{e}_3^{(0)},
\label{eq.bccW}
\end{multline}
with integer $i$, $j$ and $k$.

A $W$ point is also a point of maximum free range in the BCC lattice. The common distance
between $(0,a/4,a/2)$ and its four nearest neighbors at $(0,0,0)$, $(-a/2,a/2,a/2)$,
$(a/2,a/2,a/2)$ and $(0,0,a)$ is $\sqrt{5}a/4\approx 0.55902a$.

Since each of the 24 $W$ points is shared by two primitive unit cells of the
SC lattice, adding the $W$ points adds 12 points to the primitive unit cell,
for a total of 14 once the 2 points already present at the previous levels
of refinement are included.
An alternative way of counting what is inside the WS cell of the BCC lattice looks as follows:
each of the 24 $W$ points is shared by 4 Wigner-Seitz cells that meet at each $W$.
Adding $W$ points therefore adds $24/4=6$ points to the Wigner-Seitz cell of the
BCC lattice for a total of 7.

After these points at $W$ positions have been added to BCC lattice,
the statistics of shells of neighbors around the $\Gamma$-points 
is gathered in Table \ref{tab.nl2Gamma}.
The statistics of shells of neighbors around the $W$-points
is obviously different: Table \ref{tab.nl2W}.
\begin{table}
\caption{Statistics of distances for points with $\Gamma$-symmetry after the $W$ points
have been added to the BCC lattice.
}
\begin{tabular}{r|r}
\hline
count & squared distance \\
\hline
24 & $\frac{5}{16} a^2$ \\
8 & $\frac34 a^2$  \\
24 & $\frac{13}{16}a^2$ \\
6 & $a^2$\\
48 & $\frac{21}{16}a^2$\\
72 & $\frac{29}{16}a^2$ \\
\hline
\end{tabular}
\label{tab.nl2Gamma}
\end{table}

\begin{table}
\caption{Statistics of distances for points with $W$-symmetry after the $W$ points
have been added to the BCC lattice.
}
\begin{tabular}{r|r}
\hline
count & squared distance \\
\hline
4 & $\frac18 a^2$ \\
2 & $\frac14 a^2$ \\
4 & $\frac{5}{16}a^2$ \\
8 & $\frac38 a^2$ \\
4 & $\frac12 a^2$ \\
8 & $\frac58 a^2$ \\
8 & $\frac34 a^2$ \\
\hline
\end{tabular}
\label{tab.nl2W}
\end{table}

The weights of the two kinds of points in numerical integration
are the volumes of their Voronoi cells,
computed in Appendix \ref{app.volW}:
\begin{equation}
V_\Gamma^{(2)} = \frac{125}{1152}a^3;\quad
V_W^{(2)} =  \frac{451}{6912}a^3 ;\quad
2V_\Gamma^{(2)} + 12 V_W^{(2)} = V^{(0)}.
\label{eq.bccWv}
\end{equation}

The generic algorithm to determine the volumes of polyhedra is to 
acquire the Cartesian coordinates of all vertices, to define the 
face set as a set of co-planar triangles and to sum the contribution
of each triangle (with outwards orientation of the face normal) to the 
volume with the divergence theorem. The contribution of each triangle
is a sixth of the scalar triple product of the three vectors 
from the origin of coordinates to the triangle's
vertices.

\subsubsection{Level 2: body-centered plus $X$}

That growth
proposed in Section \ref{sec.bcW}
for the number of points by a factor 7 compared to the BCC level
may be too drastic for some applications. So in engineering practise one could
as well add the $X$ points at mid-points of the square faces of the BCC lattice, although their
free range to their 2 nearest neighbors (the cube centers) is only $a/2$, smaller than
the free range of the $W$. The $X$ are shared
between two WS cells of the BCC lattice, so the total number of points
in the WS unit cell of the BCC lattice raises by 3 to a total of 4.

With the simplified setup proposed above, the Level 2 grid points contain
the SC points of (\ref{eq.scp}), the body-centered points of
(\ref{eq.bccp}) and---after a glance at Figure \ref{fig.WSbcc}---also the face-centered
and
edge-centered grid points.
They have distance $a/2$ to their next nearest neighbors.
In the language of Brillouin zones
one could call these points $M$ points
of the (level 0) Wigner-Seitz cell.

At that level, the full set of refined grid points is a SC lattice
with lattice constant $a/2$, volume $a^3/8$ in the unit cell,
and the same directions of its unit vectors
as the host grid of level zero. In the SC unit cell of level 0 we started
with one point per unit cell, added one point per unit cell at level 1,
added three $X$ and three $M$ points per unit cell at
level 2, such that there are now 8 lattice points per level 0 unit cell.

The benefit of that choice of adding points is equivalent to
the one illustrated two refinements of the plane in Section \ref{sec.2d}:
Adding more points at finer levels onwards is a procedure
continuing  recursively with the procedure in Section \ref{sec.bcc},
because the levels zero and three contain congruential sets of lattice points.
The analysis is in that sense complete if points at $X$ had been
added at level 2.

\subsection{Level 3: body-centered with $W$ and $\Lambda$} \label{sec.WL}
If the points at $W$ have been added in level 2, the points with maximum
minimum distance to be added in the next refinement of level 3 are
located on the $\Lambda$-line of the space diagonal at $(\frac{5}{24}a,\frac{5}{24}a,\frac{5}{24}a)$. [This is close to but not
precisely at the $L$ which is at $(\frac14,\frac14,\frac14)$.]
There are three replicas of this point with the same point group
symmetry along the space
diagonal at $(\lambda a,\lambda a,\lambda a)$ with
$\lambda=1/2-5/24=7/24$, $1/2+5/24=17/24$ and $1-5/24=19/24$.
Each point of these quartets adds equivalent positions on the other
three space diagonals; so there are $4\times 4=16$ new lattice
points added here in level 3 for a combined
total of 30 in the primitive SC unit cell.

\begin{figure}
\includegraphics[width=0.8\columnwidth]{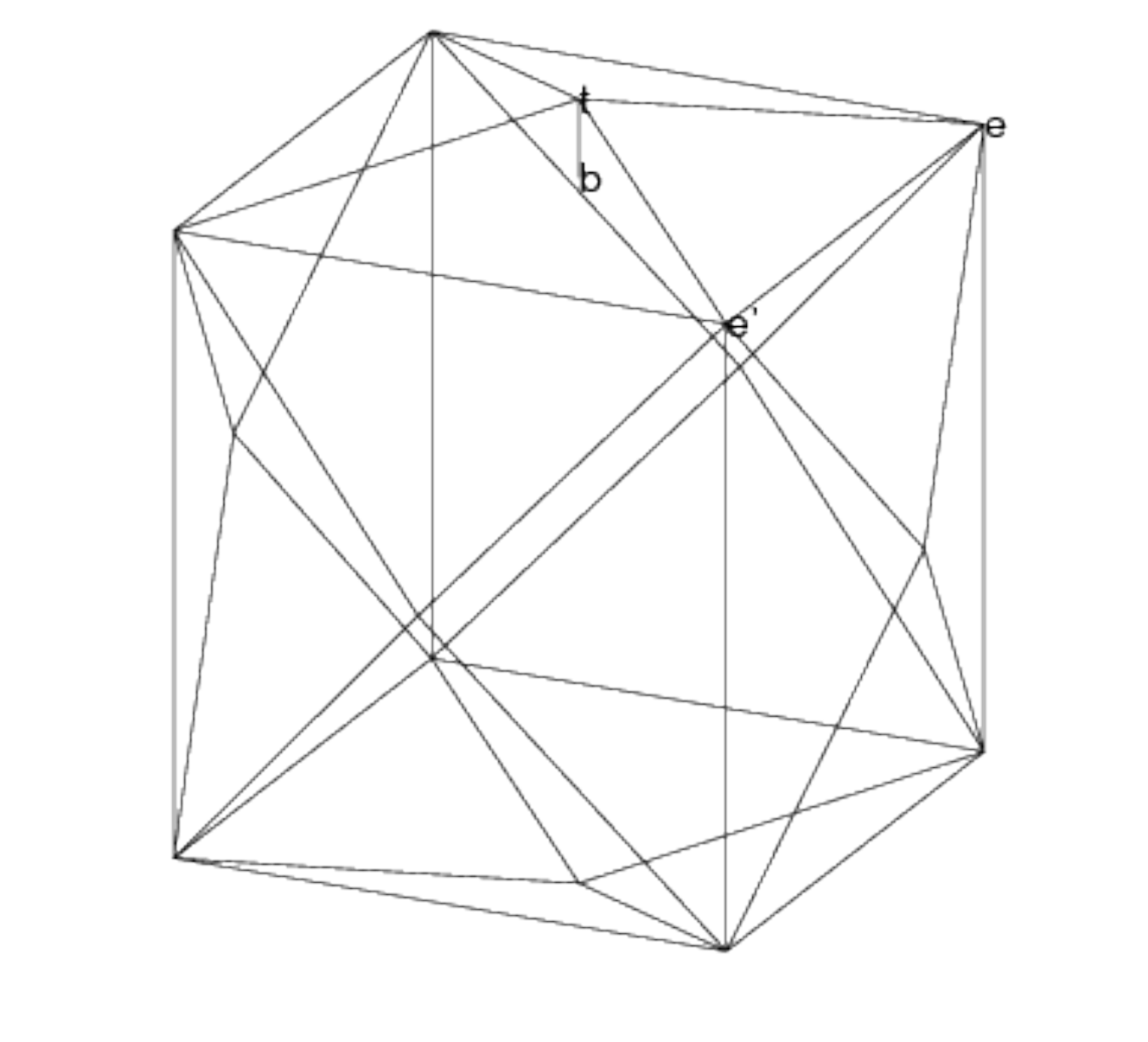}
\caption{
Geometry of the Tetrakis Hexahedra encapsulating the $\Gamma$ points of the level 2
lattice.
The edge length $e-e'$ of the cube is $5a/12$.
}
\label{fig.tOct}
\end{figure}
\begin{figure}
\includegraphics[width=0.6\columnwidth]{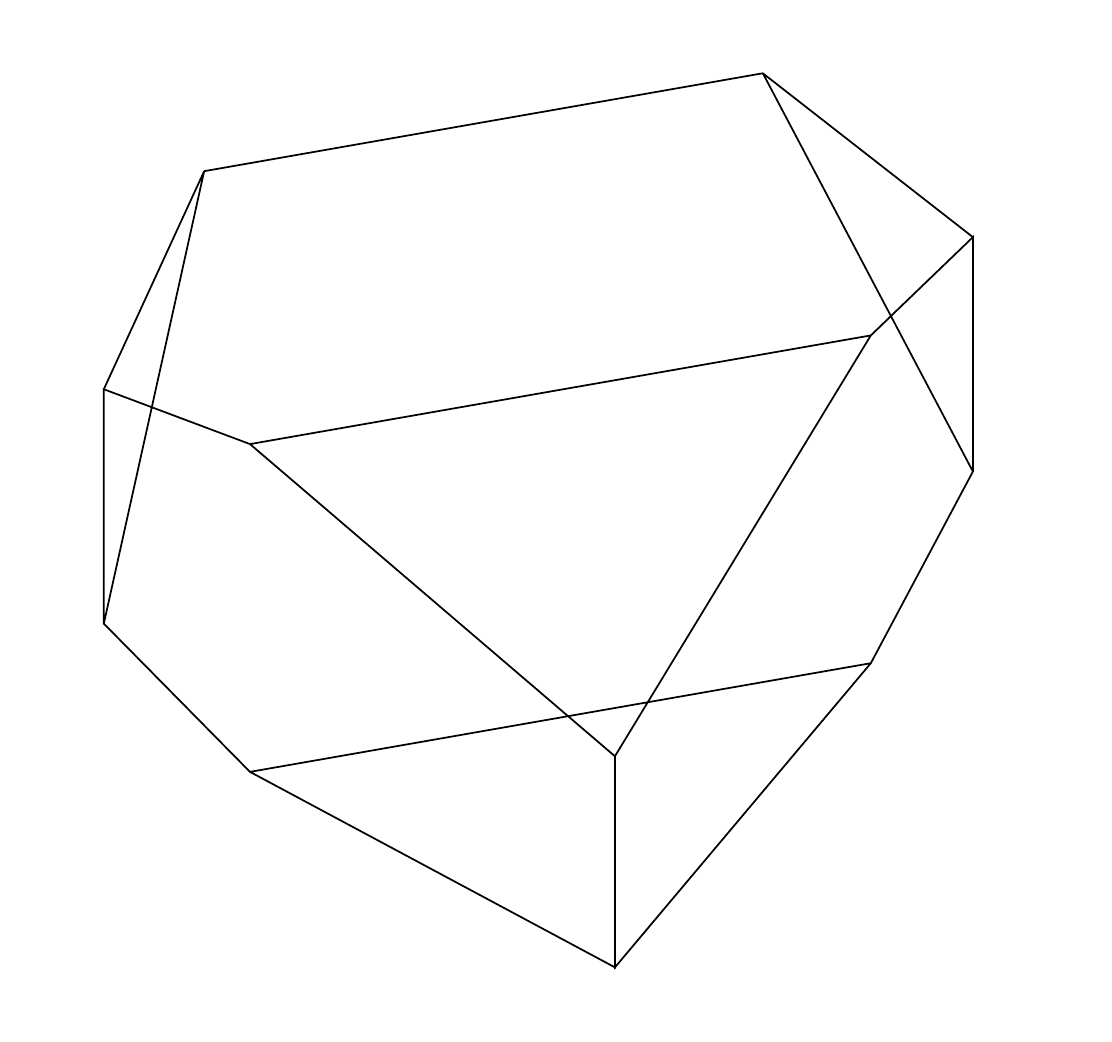}
\caption{Voronoi Octahedron encapsulating the $W$ points
of the level 2 lattice. The set of faces contains 4 regular triangles
and 4 hexahedra. Each face normal of a triangle center is a 3-fold rotation axis.
}
\label{fig.WSbccWaroundW}
\end{figure}
\begin{figure}
\includegraphics[width=0.8\columnwidth]{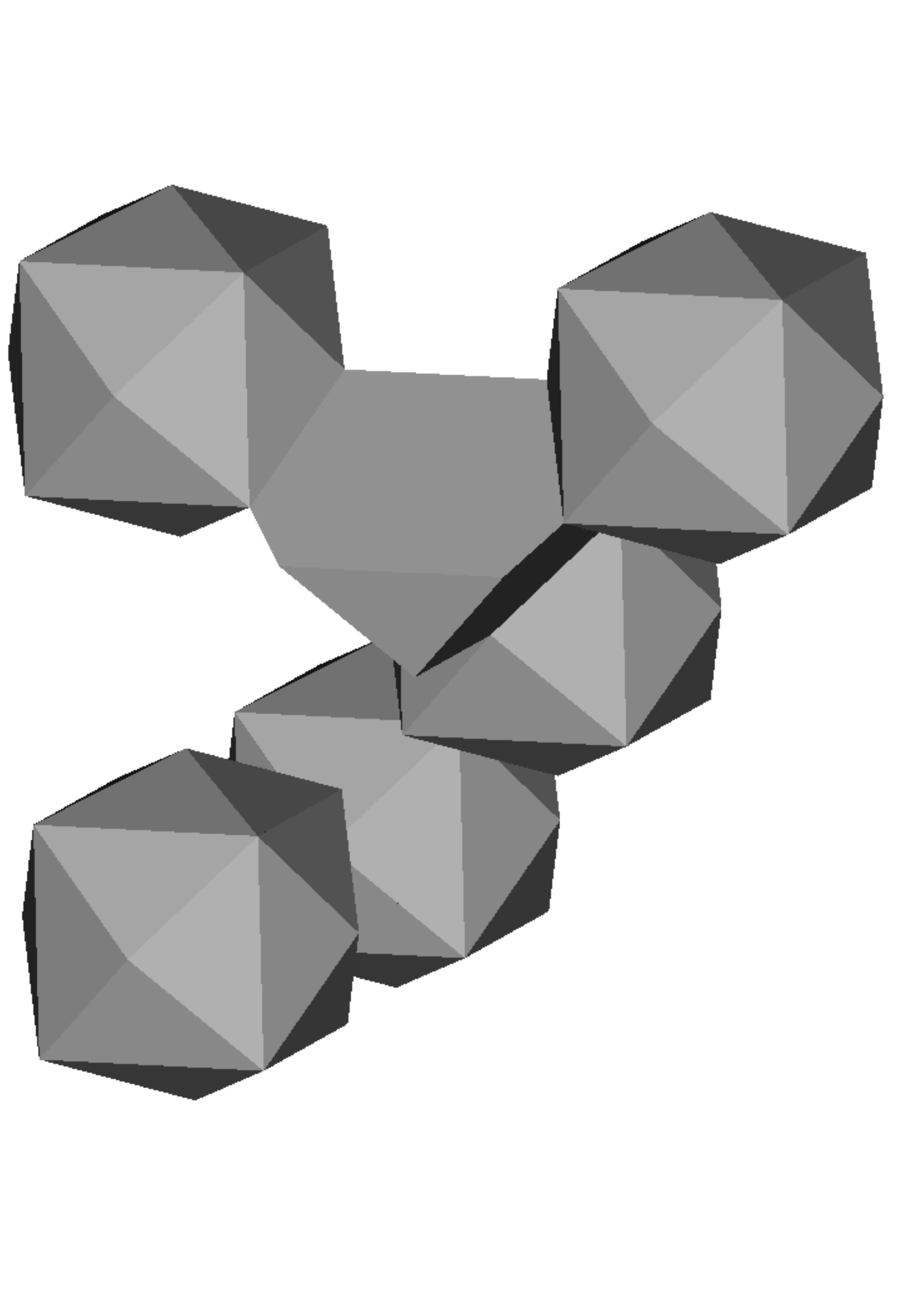}
\caption{
Level 2 dissection of the unit cell with lattice points at $\Gamma$
(the BCC lattice points) and at $W$.
The Voronoi cells around the $\Gamma$ points have the shape of Figure \ref{fig.tOct};
five of them are shown. The cells around the $W$ points have the shape of Figure \ref{fig.WSbccWaroundW};
only one of these is shown.
Each pair of cells around $\Gamma$ at a mutual distance $a$ are connected by 4
cells around the four intermediate $W$, which share the 4-fold rotation
axis that runs through $X$.
}
\label{fig.WSbccW}
\end{figure}
\begin{figure}
\includegraphics[width=0.8\columnwidth]{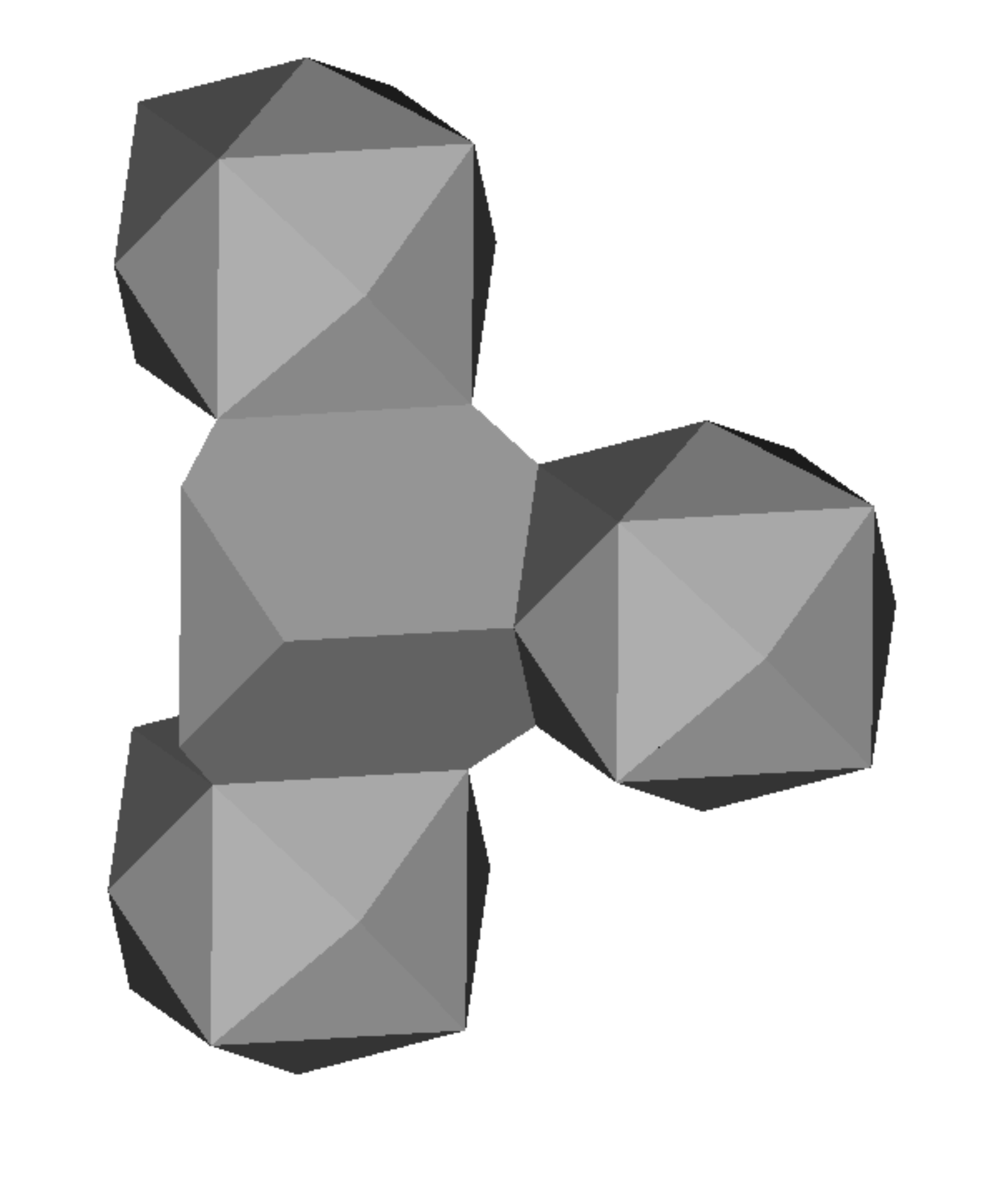}
\caption{
The Voronoi cells of Figure \ref{fig.WSbccW} from a
perspective
rotated by 90 degrees.
}
\label{fig.WSbccWalt2}
\end{figure}

The distance of this point $(5a/24,5a/24,5a/24)$ to its 7 nearest
neighbors
is $5\sqrt{3}/24a\approx 0.3608439a$.
The set of 7 nearest neighbors contains $(0,0,0)$ and 6 $W$-points of
the list (\ref{eq.bccW}),
$(0,a/4,a/2)$,
$(0,a/2,a/4)$,
$(a/4,a/2,0)$,
$(a/2,a/4,0)$,
$(a/4,0,a/2)$,
$(a/2,0,a/4)$.
The geometric interpretation of these coordinates looks as follows:
After the $W$-points have been added to the lattice,
the Voronoi cell around the $\Gamma$-points are 24-hedrons
(Tetrakis Hexahedra)
with flat square pyramids glued to each side of a cube
around the $\Gamma$-point (Figure \ref{fig.tOct}).

The cube inside has an edge length of $5a/12$, so its eight vertices
have coordinates of $\pm 5a/24$.

The Voronoi cells around the $W$-points at level 2
are 8-hedrons (Figure \ref{fig.WSbccWaroundW})
with
4 triangles that touch the triangles
of the truncated octahedra and 4 irregular hexagons that touch 8-hedrons
of adjacent $W$-points.
The short sides of the hexagons are the residual left from the distance
of $\sqrt{3} a/2$ between two truncated octahedra along the space diagonal
after removal of the diagonals inside the two cubes of the truncated
octahedra, $\sqrt{3} a/2-2\times \sqrt{3} 5a/24 = \sqrt{3}a/12\approx 0.1443376 a$.
Two such edges are the two rightmost edges of the 8-hedron
in Figure \ref{fig.WSbccWalt2} that each join two of the Tetrakis Hexahedra.

The points added here at level 3
are those vertices around the $W$ Voronoi cell
where a triangle meets a \emph{short} edge of
the hexagon. At these points the Tetrakis Hexahedron of the
Voronoi cell around the $\Gamma$-points meets 6 Voronoi cells
of $W$ points. They lie on the axes with 3m symmetry (along the
space diagonal) of the SC lattice.

The factor $5/24$ that fixes the position of
the new $\Lambda$ points is then found by noticing that these
short edges run along a space diagonal, that the projection
of the new point on the $xy$ plane runs along the plane diagonal,
and to solve the planar Voronoi cell problem for
the triple of points $(0,0)$, $(a/2,a/4)$ and $(a/4,a/2)$
in that plane.

The Voronoi cell around a $\Gamma$ point at level 3
is an Octahedron with
space diagonal $5a/8$ (all bounding box coordinates are $\pm 5a/16$),
so the edge length at the square base of the pyramid is $5a/(8\sqrt 2)$,
and the volume is \cite{ConwayITIT28}
\begin{equation}
V_\Gamma^{(3)}=\frac{125}{3072}a^3\approx 0.04069010 a^3.
\label{eq.gamma3v}
\end{equation}

The degeneracies of the various points in the SC lattice require
\begin{equation}
2V_\Gamma^{(3)}+16V_\Lambda^{(3)}+12V_W^{(3)}=a^3.
\end{equation}
Taking the volume $V_\Lambda^{(3)}$ from Eq. \ (\ref{eq.Vlam3})
and solving for $V_W^{(3)}$ we obtain
\begin{equation}
V_W^{(3)} = \frac{24505}{663552}a^3\approx 0.03693003a^3
\label{eq.W3v}
\end{equation}
for the volumes associated with the $W$ points at level 3.

\begin{figure}
\includegraphics[width=0.95\columnwidth]{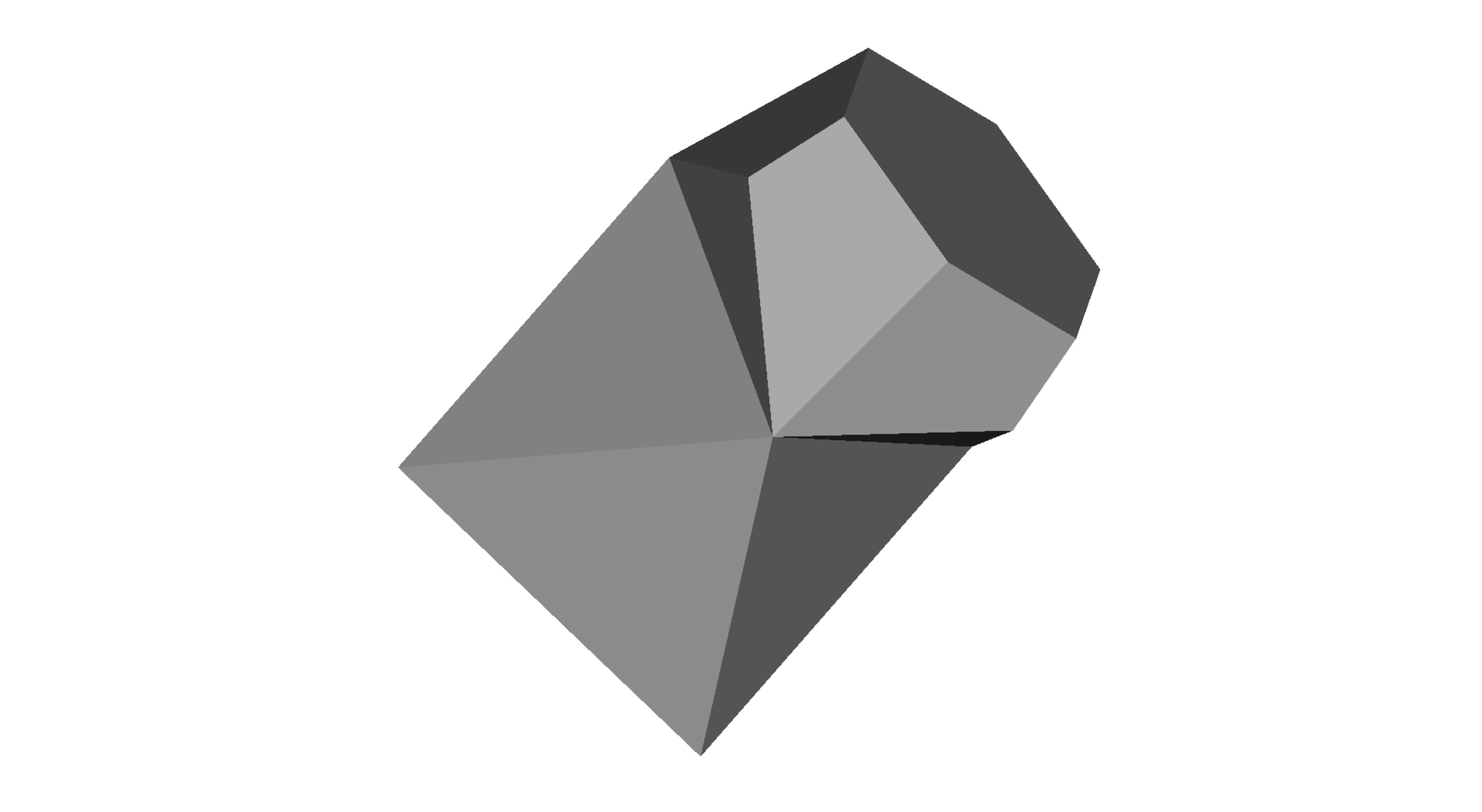}
\caption{Level 3 Voronoi cells around $\Gamma$ (Octahedron) and around 
a neighboring $\Lambda$. The $\Lambda$ cell has 11 faces: 1 regular hexagonal
side where it points along the space diagonal (shared with another $\Lambda$ cell),
1 regular triangular side
where it attaches to the Octahedron, and 6 quadrangles (shared with $W$ cells)
and 3 triangles (shared with $\Lambda$ cells)
at the interface. The symmetry is the 3-fold rotation axis along the space
diagonal.
}
\label{fig.gammWLlam1}
\end{figure}

\begin{figure}
\includegraphics[width=0.8\columnwidth]{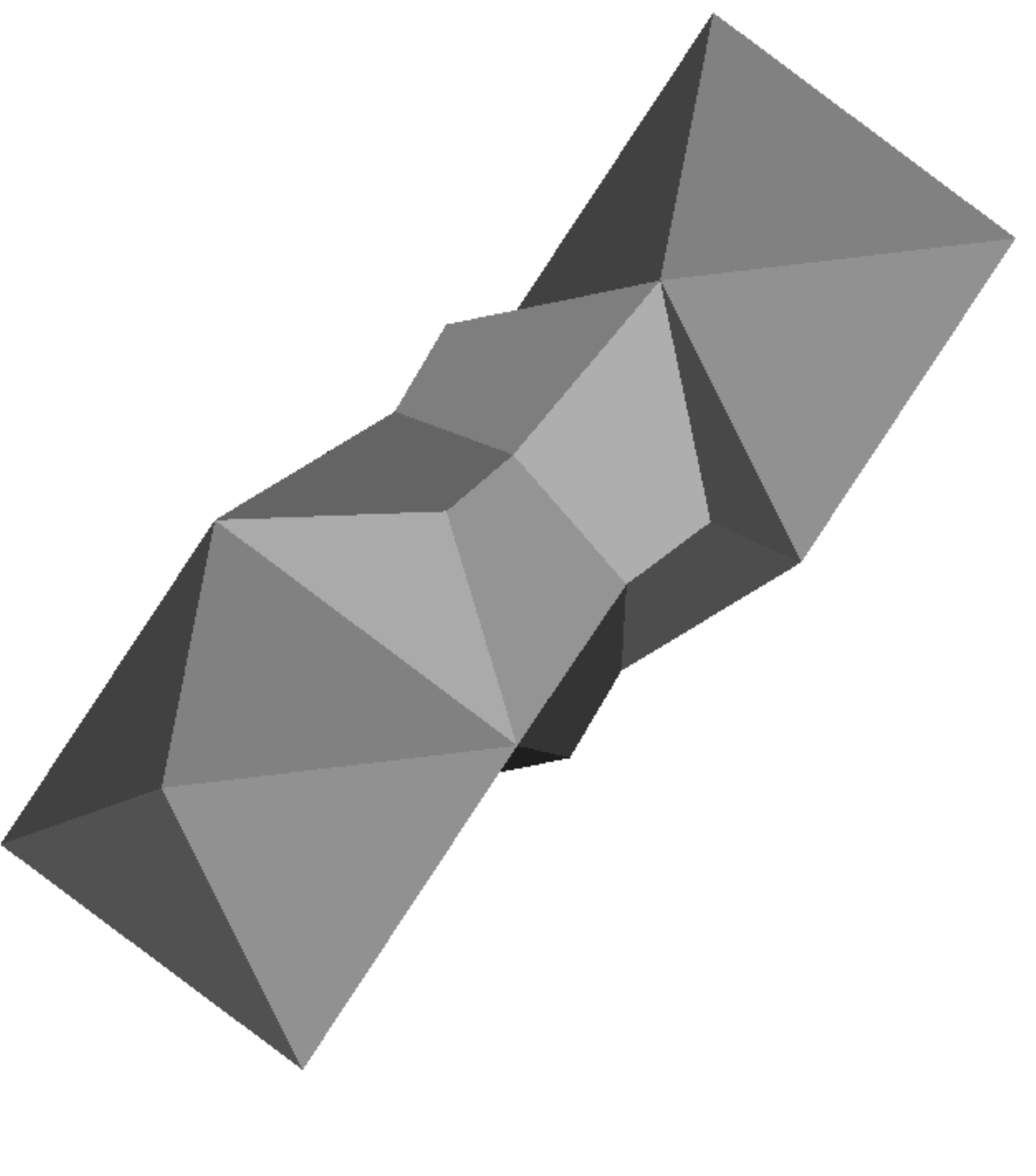}
\caption{Illustration of the bridging between the Octahedra of the
level 3 lattice along each space diagonal of the SC lattice by pairs
of Voronoi cells around the $\Lambda$ points. This follows by extension
of Figure \ref{fig.gammWLlam1} with mirror symmetry at the hexagonal face
of the $\Lambda$ cell.
}
\label{fig.gammWLcross}
\end{figure}

\begin{figure}
\includegraphics[width=0.4\columnwidth]{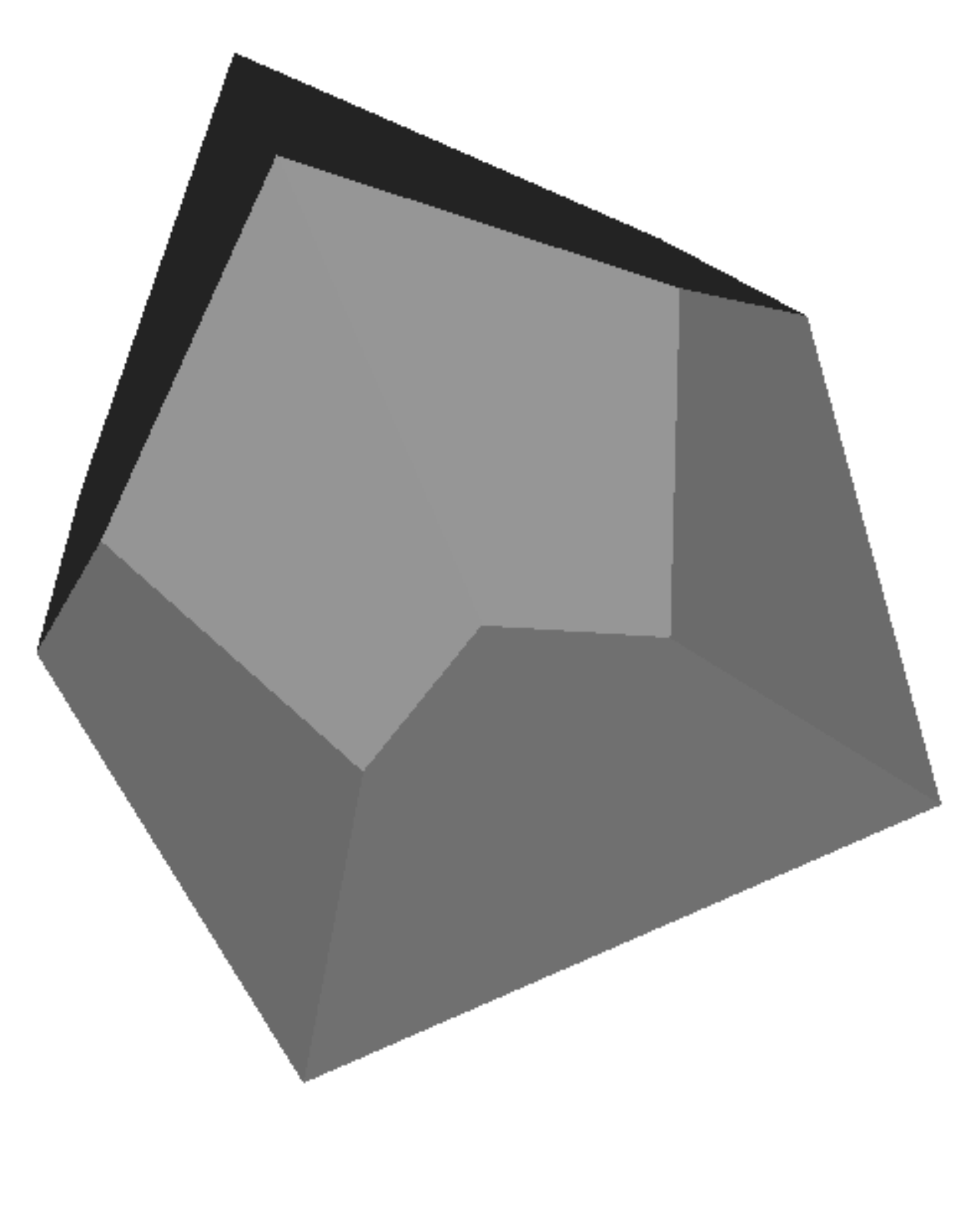}
\includegraphics[width=0.4\columnwidth]{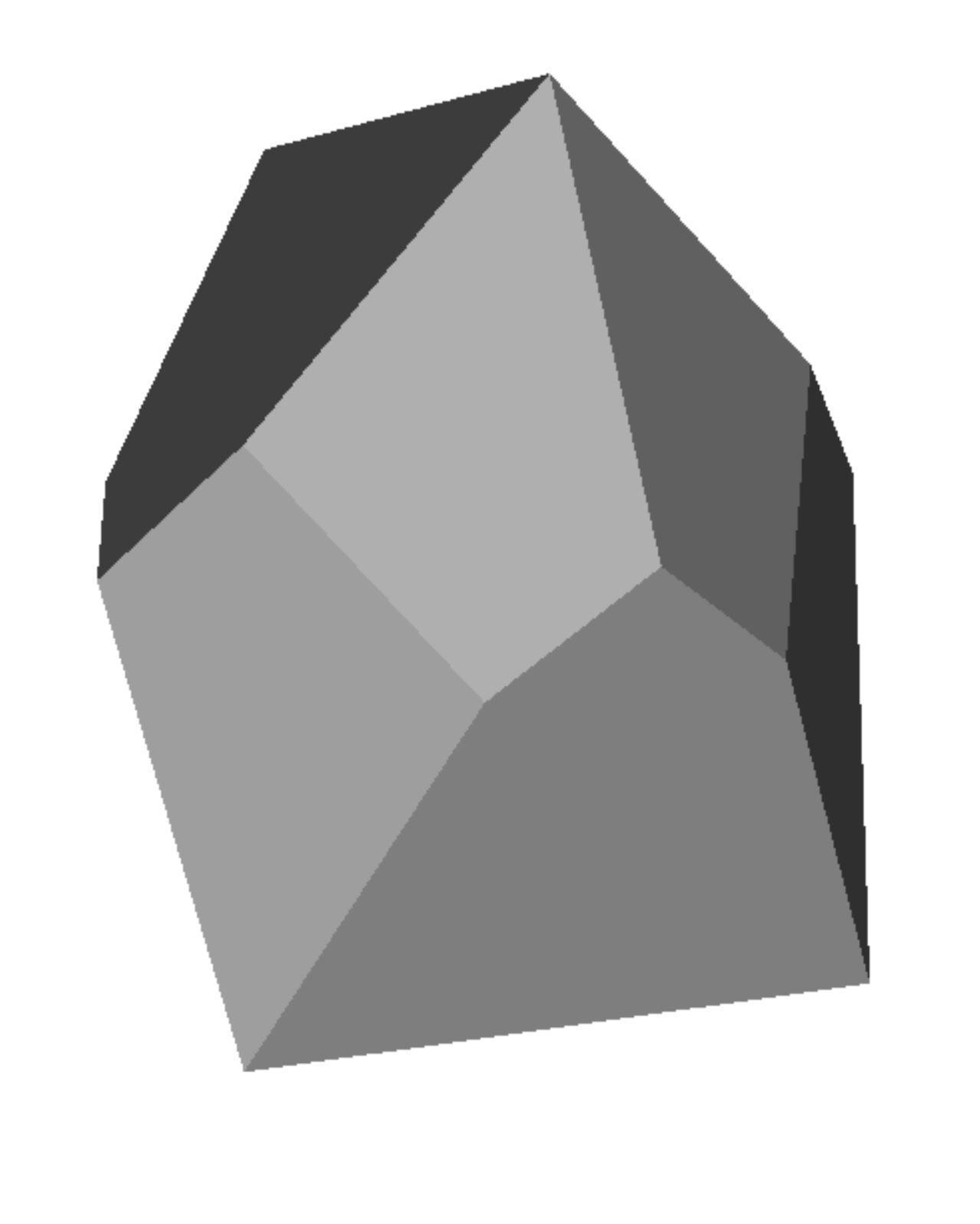}
\caption{Level 3 Voronoi cell around any of the $W$ points after $W$ and $\Lambda$
points have been added to the BCC lattice.
It has 12 faces:
two pairs of pentagons, each pair with a common long edge (right dihedral
angle), plus a band
of eight quadrangles
that runs in between.
}
\label{fig.gammWLWup}
\end{figure}

\begin{figure}
\includegraphics[width=0.8\columnwidth]{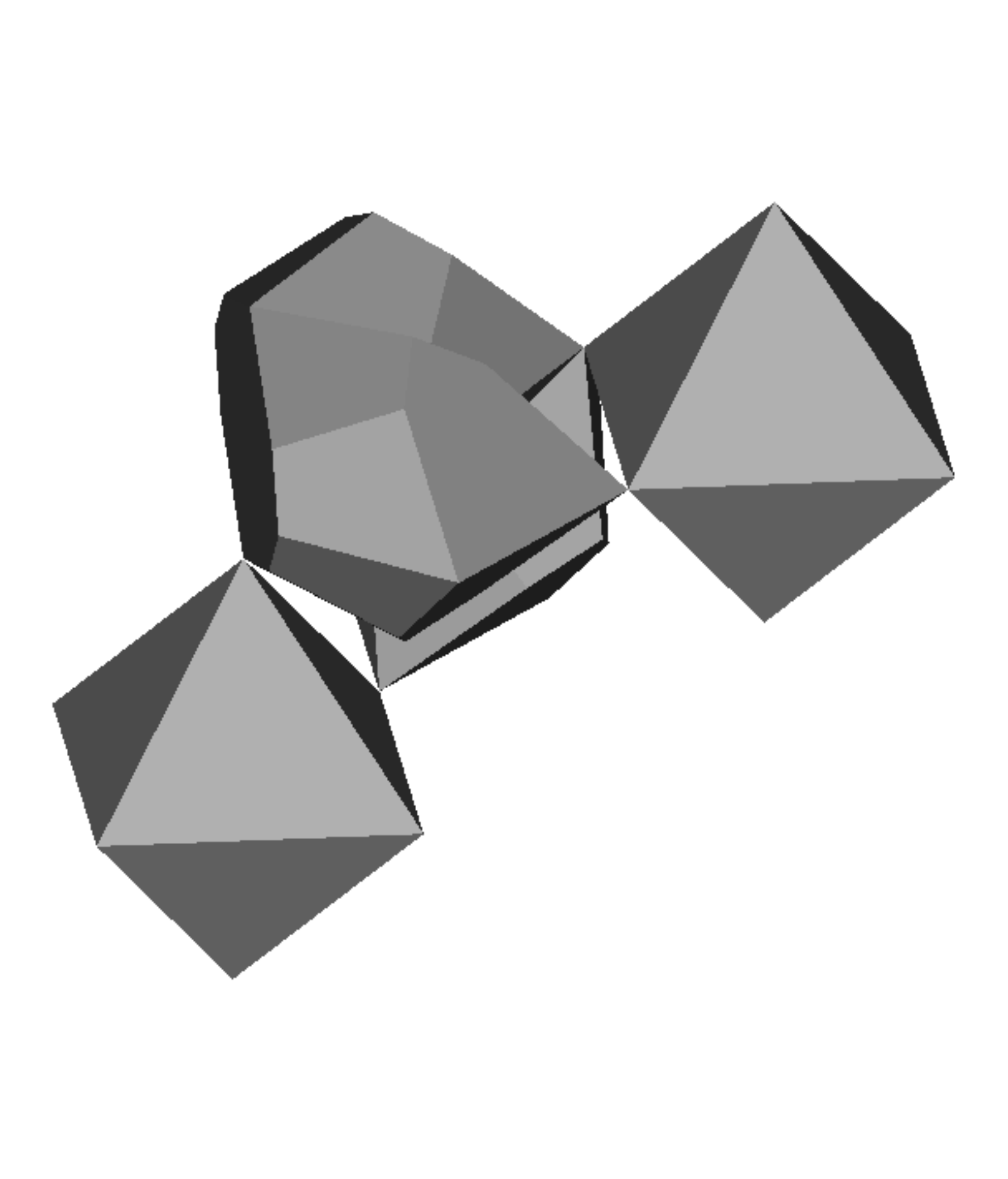}
\caption{Voronoi cells at a two $\Gamma$ points plus three $W$ cells 
of the shape of Figure \ref{fig.gammWLWup} around
$W$ points after $W$ and $\Lambda$
points have been added to the BCC lattice.
Each of the long edges of the pentagons of Figure \ref{fig.gammWLWup} touches
one point of an Octahedron.
}
\label{fig.gammWLW3}
\end{figure}

\begin{figure}
\includegraphics[width=0.8\columnwidth]{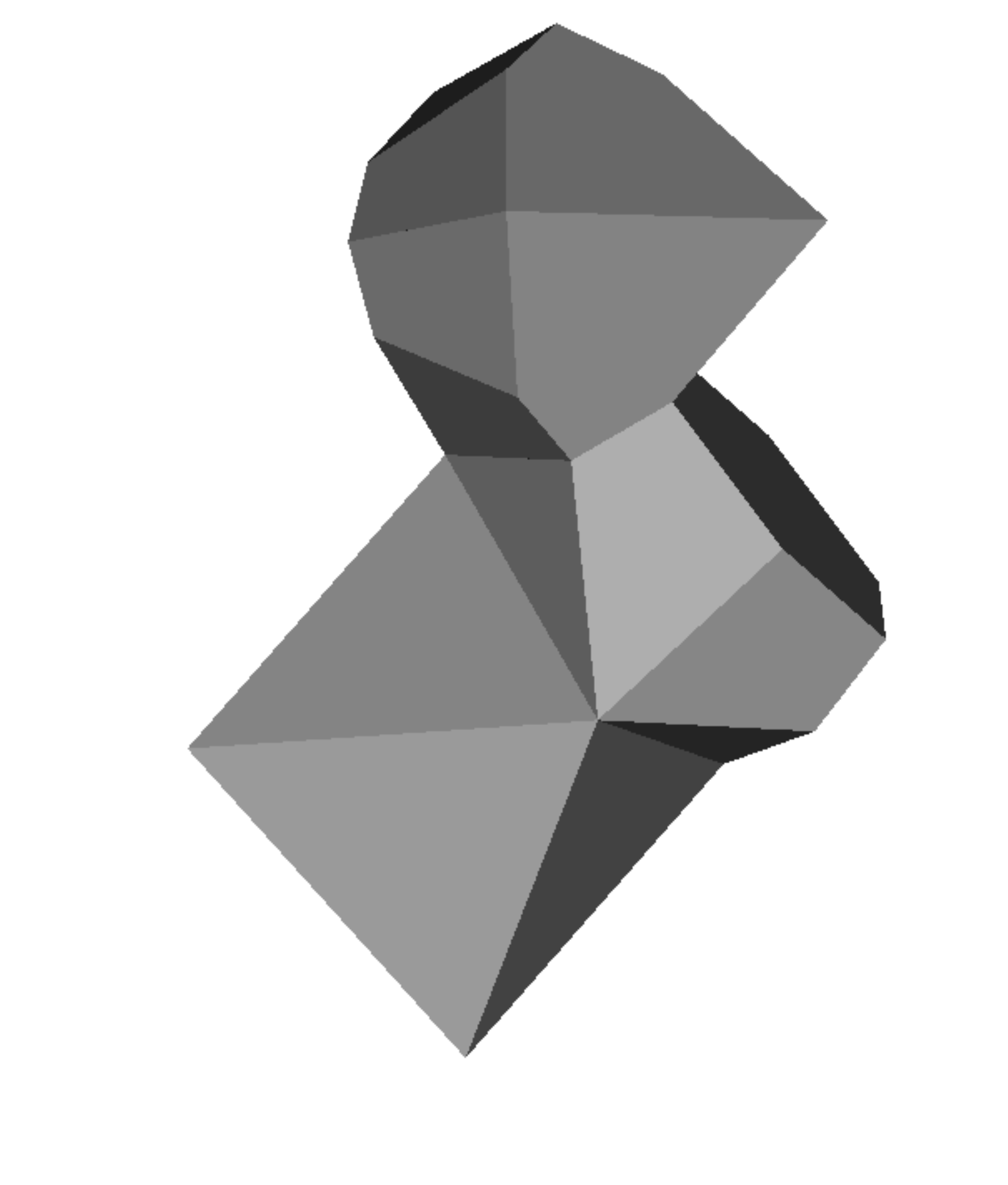}
\caption{The configuration if a single cell of Figure \ref{fig.gammWLWup}
is added to Figure \ref{fig.gammWLlam1}.
}
\label{fig.gammWLlam1Wup}
\end{figure}

\begin{figure}
\includegraphics[width=0.8\columnwidth]{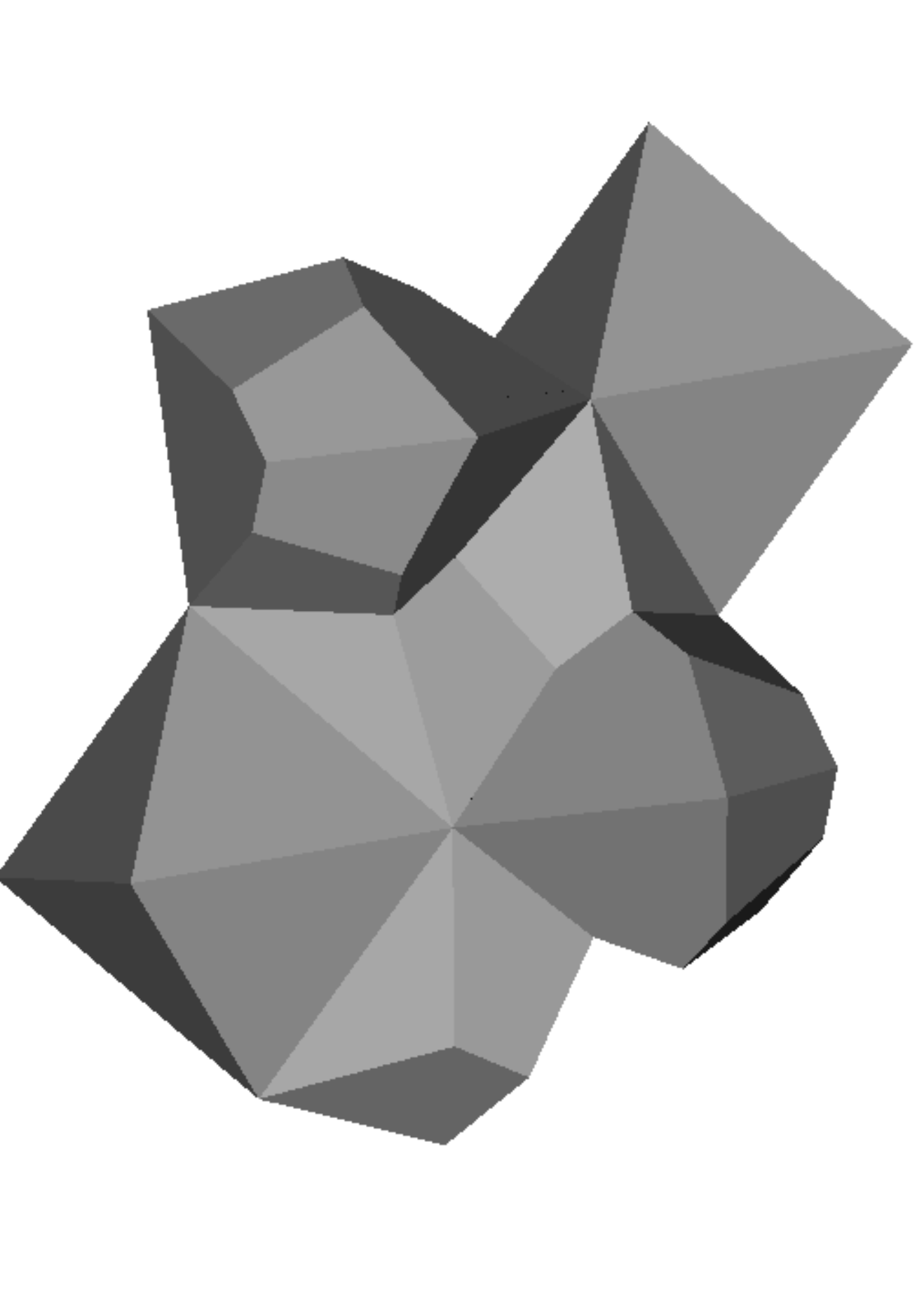}
\caption{
A composite of Figure \ref{fig.gammWLcross}, Figure \ref{fig.gammWLlam1Wup}
and two further cells
(around a $\Lambda$ point on a different space diagonal and around
a second $W$ point).
}
\label{fig.gammWLall}
\end{figure}

Finally, Figures \ref{fig.gammWLlam1}--\ref{fig.gammWLall}
try to give an impression of how the three
different types of Voronoi cells at level 3 subdivide the
region near two points of the BCC WS cell.

\section{Summary}
Regular dissection of space with Voronoi cells on the points of
the simple cubic grid assigns the volume (\ref{eq.scv}) to each cell.
If points are added at the body centers (\ref{eq.bccp}),
the volume (\ref{eq.bccv}) is assigned to the Voronoi cells around the grid points.
If another set of points is added at all $W$ positions (\ref{eq.bccW}),
Voronoi cells of two different shapes appear, characterized by 
volumes defined in Eq.\ (\ref{eq.bccWv}).
If finally that set is augmented by points on the $\Lambda$ line,
space is divided by three different types of Voronoi cells with volumes
of (\ref{eq.Vlam3}), (\ref{eq.gamma3v}) and (\ref{eq.W3v}), respectively.

\begin{acknowledgments}
Many images of Voronoi cells shown here have been created by loading
stereolithography data into \texttt{MeshLab} \cite{meshlab}.
\end{acknowledgments}

\appendix
\section{Shape of the Voronoi cells at level 2} \label{app.volW}
\subsection{Tetrakis Hexahedra at $\Gamma$} \label{sec.tetra}

The height of the six pyramids  that cover the Tetrakis hexahedron in Figure \ref{fig.tOct}
is determined by considering the mutual tilt of two planes that 
dissect the space from a $\Gamma$-point towards to neighboring $W$-points,
for example $\Gamma$ at $(a/2,a/2,a/2)$, $W$ at $(0,a/2,a/4)$ and
$W$ at $(0,a/2,3a/4)$.
The distances $\Gamma-W$ are 
$\sqrt{(1/2)^2+0^2+(1/4)^2}a = \sqrt{5}a/4$.
The angle of the line $\Gamma\to W$ versus the horizontal is 
$\tan \phi = (a/4)/(a/2)$, $\phi = \arctan\frac{1}{2} \approx 26.5650^\circ \approx 0.463647$ rad
\cite[A073000]{EIS}.
This is the inclination of the triangles of the pyramids of the Tetrakis
Hexahedron versus their cubic base, the tilt of the triangle $e'-e-t$
versus the horizontal in Figure \ref{fig.tOct}.
The distance from the base $b$ to the top $t$ of the pyramids is
a quarter of the edge length, $5a/48$, because the slope of $\phi$ is $1/2$. The distance
$e-t$ is $\sqrt{(5/24)^2+(5/24)^2+(5/48)^2}a=5a/16=0.3125a$.
This yields a height of the planar triangle of 
$\sqrt{(5/16)^2-(5/24)^2}a=5^{3/2}a/48\approx 0.23292a$.
The internal angles of the triangle are
$\arccos\frac{5a/24}{5a/16}=\arccos\frac23\approx 0.8411$ rad
$\approx 48.189685^\circ$ for $e-e'-t$ \cite[A228496]{EIS}
and the 180$^\circ$ complement $83.62063^\circ$ for $e'-t-e$.
The volume of each pyramid is a third of the product of height by base area,
$(1/3)\times (5a/48)\times (5a/12)^2=
125a^3/20736\approx 0.00603a^3$.
The 6 pyramids plus the cube encapsulate a volume of
\begin{equation}
V^{(2)}_\Gamma = \frac{125}{1152}a^3\approx 0.1085a^3.
\label{eq.V2gamma}
\end{equation}

\subsection{Truncated Tetrahedra at $W$}

In the (irregular) truncated tetrahedra around the $W$ points in Figure
\ref{fig.WSbccWaroundW},
the edge lengths of the triangular surfaces are the edge lengths
$5a/12$ and $5a/16$ of the pyramids determined in Section \ref{sec.tetra}.
They also fix edge lengths $5a/12$ and $5a/16$
of the four hexagons.
The missing edge length is defined by subtracting two pyramidal
heights and two cube half edges
from the distance $a$ between two $\Gamma$-points of the
lattice, proposed by the upper horizontal edge that connects
two pyramidal summits in Figure \ref{fig.WSbccW},
$a-2\times 5a/24- 2\times 5a/48=3a/8 =0.375a$. In summary, the
hexagonal faces of the cell around the $W$ points have
one edge of length
$5a/12$,
two edges of
length
$5a/16$,
two edges of length
$\sqrt{3}a/12$,
and
one edge of length
$3a/8$.

One can reconstruct all the hexagonal surfaces starting with
the vertex set of one of them,  for example
$(5/24)(\mathbf{e}_1^{(0)}
+\mathbf{e}_2^{(0)}
+\mathbf{e}_3^{(0)})$,
$(7/24)(\mathbf{e}_1^{(0)}
+\mathbf{e}_2^{(0)}
+\mathbf{e}_3^{(0)})$,
$(5/24)(-\mathbf{e}_1^{(0)}
+\mathbf{e}_2^{(0)}
+\mathbf{e}_3^{(0)})$,
$(7/24)(-\mathbf{e}_1^{(0)}
+\mathbf{e}_2^{(0)}
+\mathbf{e}_3^{(0)})$,
$(3/16)\mathbf{e}_1^{(0)}
+(1/2)\mathbf{e}_2^{(0)}
+(1/2)\mathbf{e}_3^{(0)}$,
$-(3/16)\mathbf{e}_1^{(0)}
+(1/2)\mathbf{e}_2^{(0)}
+(1/2)\mathbf{e}_3^{(0)}$,
and applying the group operations of the BCC space group
to it.

The volume of the Voronoi cell is determined by
partitioning of the SC cell volume $a^3$ in (\ref{eq.scv}):
The volume of the cells around $\Gamma$ points is (\ref{eq.V2gamma}).
There are effectively two of them in the SC unit cell (one in the
BCC body center and an eights of each of the shared points
in the eight corners).
There are effectively 12 cells around $W$ per SC unit cell. Each
of them needs to contribute
\begin{equation}
V_W^{(2)}=\frac{451}{6912}a^3 \approx 0.0652488a^3
\end{equation}
to satisfy equation (\ref{eq.bccWv}).

\section{Voronoi Cells of $\Lambda$ at level 3}

The coordinates of points of the $\Lambda$ cell in Figure \ref{fig.gammWLlam1} are determined
as follows:
\begin{itemize}
\item
The points on the triangular face attached to the Octahedron at $\Gamma$ are 
fixed by the Octahedron,
$(5/16)\mathbf{e}_3^{(0)}$,
$(5/16)\mathbf{e}_1^{(0)}$,
and $(5/16)\mathbf{e}_2^{(0)}$.
\item
A first vertex on the hexagonal face is at
$(95/288)\mathbf{e}_1^{(0)}
+(13/144)\mathbf{e}_2^{(0)}
+(95/288)\mathbf{e}_3^{(0)}$.
It is found by construction of (i) the plane than divides
the space between $\Lambda$ at
$(5/24)(\mathbf{e}_1^{(0)}
+\mathbf{e}_2^{(0)}
+\mathbf{e}_3^{(0)}
)$
and $\Lambda$ at
$(7/24)(\mathbf{e}_1^{(0)}
+\mathbf{e}_2^{(0)}
+\mathbf{e}_3^{(0)}
)$, of (ii) the plane that divides
the space between $\Lambda$ at
$(5/24)(\mathbf{e}_1^{(0)}
+\mathbf{e}_2^{(0)}
+\mathbf{e}_3^{(0)}
)$
and $W$ at
$(1/4)\mathbf{e}_1^{(0)}
+(1/2)\mathbf{e}_3^{(0)}$, and (iii) of the plane that divides
the space between $\Lambda$ at
$(5/24)(\mathbf{e}_1^{(0)}
+\mathbf{e}_2^{(0)}
+\mathbf{e}_3^{(0)}
)$
and $W$ at
$(1/2)\mathbf{e}_1^{(0)}
+(1/4)\mathbf{e}_3^{(0)}$, followed by intersection of these three planes.
The other five points are obtained from this one by rotations around the space
diagonal by multiples of $60^\circ$,
$(59/144)\mathbf{e}_1^{(0)}
+(49/288)\mathbf{e}_2^{(0)}
+(49/288)\mathbf{e}_3^{(0)}$,
$(95/288)\mathbf{e}_1^{(0)}
+(95/288)\mathbf{e}_2^{(0)}
+(13/144)\mathbf{e}_3^{(0)}$,
$(49/288)\mathbf{e}_1^{(0)}
+(59/144)\mathbf{e}_2^{(0)}
+(49/288)\mathbf{e}_3^{(0)}$, and so on.

\item
The missing three points that are vertices of quadrangles are at
$(35/128)(\mathbf{e}_1^{(0)}
+\mathbf{e}_3^{(0)})$,
$(35/128)(\mathbf{e}_1^{(0)}
+\mathbf{e}_2^{(0)})$, and
$(35/128)(\mathbf{e}_2^{(0)}
+\mathbf{e}_2^{(0)})$.
The quadrangle is constructed by (i) the plane than divides
the space between $\Lambda$ at
$(5/24)(\mathbf{e}_1^{(0)}
+\mathbf{e}_2^{(0)}
+\mathbf{e}_3^{(0)}
)$
and $\Lambda$ at
$(5/24)(\mathbf{e}_1^{(0)}
-\mathbf{e}_2^{(0)}
+\mathbf{e}_3^{(0)}
)$, (ii) the plane that divides
the space between $\Lambda$ at
$(5/24)(\mathbf{e}_1^{(0)}
+\mathbf{e}_2^{(0)}
+\mathbf{e}_3^{(0)}
)$
and $W$ at
$(1/4)\mathbf{e}_1^{(0)}
+(1/2)\mathbf{e}_3^{(0)}$, (iii) the plane that divides
the space between $\Lambda$ at
$(5/24)(\mathbf{e}_1^{(0)}
+\mathbf{e}_2^{(0)}
+\mathbf{e}_3^{(0)}
)$
and $W$ at
$(1/2)\mathbf{e}_1^{(0)}
+(1/4)\mathbf{e}_3^{(0)}$, followed by intersection of these three planes.
\end{itemize}

The volume of the Voronoi cell surrounding
$\Lambda$ points
is computed by the divergence theorem, adding up the 
contributions of the eleven faces
\begin{itemize}
\item
one regular hexagon,
\item
six quadrangles that share
sides with the hexagon,
\item
one triangle shared with a face of the octahedron,
\item
and three triangles that share one side with the octahedron and two sides
with two of the quadrangles.
\end{itemize}
The result is
\begin{equation}
V_\Lambda^{(3)}=\frac{26291}{884736}a^3\approx
0.0297162a^3.
\label{eq.Vlam3}
\end{equation}

\bibliographystyle{amsplain}
\bibliography{all}

\providecommand{\bysame}{\leavevmode\hbox to3em{\hrulefill}\thinspace}
\providecommand{\MR}{\relax\ifhmode\unskip\space\fi MR }
\providecommand{\MRhref}[2]{%
  \href{http://www.ams.org/mathscinet-getitem?mr=#1}{#2}
}
\providecommand{\href}[2]{#2}
\begin{thebibliography}{10}

\bibitem{Bouckaert}
L.~P. Bouckaert, R.~Smoluchowski, and E.~Wigner, \emph{Theory of {Brillouin}
  zones and symmetry properties of wave functions in crystals}, Phys.\ Rev.
  \textbf{50} (1936), no.~1, 58--67.

\bibitem{ConwayITIT28}
John~H. Conway and Neil J.~A. Sloane, \emph{Voronoi regions of lattices, second
  moments of polytopes, and quantization}, IEEE Trans. Inf. Theory \textbf{28}
  (1982), no.~2, 211--226.

\bibitem{CunninghamPRB10}
S.~L. Cunningham, \emph{Special points in the two-dimensional {Brillouin}
  zone}, Phys. Rev. B \textbf{10} (1974), no.~12, 4988--4994.

\bibitem{ElliottPR96}
R.~J. Elliott, \emph{Spin-orbit coupling in band theory---character tables for
  some ``double'' space groups}, Phys. Rev. \textbf{96} (1954), 280.

\bibitem{EntezariIEEEvis04}
Alireza Entezari, Ramsay Dyer, and Torsten M\"oller, \emph{Linear and cubic box
  splines for the body centered cubic lattice}, IEEE Visualization 2004, IEEE,
  2004, pp.~11--18.

\bibitem{HerringLN}
C.~Herring, \emph{Accidental degeneracy in the energy bands of crystals},
  Lecture Notes and Supplements in Physics, vol.~5, W. A. Benjamin, New York,
  Amsterdam, 1964, pp.~240--248.

\bibitem{HerringPR52}
Conyers Herring, \emph{Accidental degeneracy in the energy bands of crystals},
  Phys.\ Rev. \textbf{52} (1937), no.~4, 365--373, Reprint: \cite{HerringLN}.

\bibitem{meshlab}
{ISTI Visual Computing Lab}, \emph{Meshlab, a system for processing and editing
  of unstructured triangular meshes}, 2012.

\bibitem{PackPRB16}
James~D. Pack and Hendrik~J. Monkhorst, \emph{``special points for
  {Brillouin}-zone integrations''.---a reply}, Phys.\ Rev.\ B \textbf{16}
  (1977), no.~4, 1748--1749.

\bibitem{EIS}
Neil J.~A. Sloane, \emph{The {O}n-{L}ine {E}ncyclopedia {O}f {I}nteger
  {S}equences}, Notices Am.\ Math.\ Soc. \textbf{50} (2003), no.~8, 912--915,
  http://oeis.org/. \MR{1992789 (2004f:11151)}

\bibitem{TerzibaschPSS133}
T.~Terzibaschian and R.~Enderlein, \emph{The irreducible representations of the
  two-dimensional space group of crystal surfaces}, Phys. Status Solidi B
  \textbf{133} (1986), no.~2, 443--461.

\end{thebibliography}

\end{document}